\documentclass{amsart}

\usepackage{amsmath}
\usepackage{amssymb}
\usepackage{amsopn}
\usepackage{amscd}
\usepackage{epsfig}
\usepackage{amsfonts}
\usepackage{latexsym}
\usepackage{graphicx}
\usepackage{color}
\usepackage{tikz}
\usepackage{wrapfig}
\usetikzlibrary{calc}

\newtheorem{theorem}{Theorem}

\newtheorem{proposition}{Proposition}
\newtheorem{lemma}{Lemma}

\definecolor{trp}{rgb}{1,1,1}

\definecolor{red}{rgb}{1,0,.2}

\definecolor{blue}{rgb}{0,0,1}

\definecolor{rgrey}{rgb}{.8,0.4,.4}  
\definecolor{grey}{rgb}{.13,.13,.13}  

\definecolor{green}{rgb}{0.0,0.4,0.2}

\newcommand{\MM}[4]{$\left( \begin{smallmatrix} #1 & #2\\[.05cm] #3 & #4 \end{smallmatrix} \right)$ }

\newcommand{\pq}[2]{ \left( \begin{smallmatrix} #1\\ #2 \end{smallmatrix} \right)}
\newcommand{\PQ}[2]{ \left( \begin{matrix} #1\\ #2 \end{matrix} \right)}

\newcommand{\mor}[2]{\Big\{
   \begin{aligned} 0 & \rightarrow #1\\[-.1cm] 1 & \rightarrow #2  \end{aligned}}

\title[]{Substitution invariant Sturmian words and binary trees}
\author[]{Michel Dekking}

\date{\today}

\begin{document}

\maketitle

\smallskip

\begin{center}
{ \footnotesize \noindent Delft University of Technology, \\Faculty EEMCS, P.O.~Box 5031,\\ 2600 GA
Delft, The Netherlands.\\ \;Email: F.M.Dekking@math.tudelft.nl }
\end{center}


\begin{abstract}
We take a global view at substitution invariant Sturmian sequences. We show that homogeneous substitution invariant Sturmian sequences $s_{\alpha,\alpha}$ can be indexed by two binary trees, associated directly to Johannes Kepler's tree of harmonic fractions from 1619. We obtain similar results for the inhomogeneous sequences $s_{\alpha,1-\alpha}$ and $s_{\alpha,0}$.

\medskip  

\noindent{\em Key words.}{\;Sturmian word; morphism; substitution; Sturm numbers; binary tree; harmonic fraction; Kepler.} 

\medskip


\end{abstract}

\bigskip

\section{Introduction}

A Sturmian word $w$ is an infinite word $w=w_1w_2\dots$, in which occur only $n+1$ subwords of length $n$ for $n=0,1,2\dots$. It is well known (see, e.g.~\cite{Lothaire}) that the Sturmian words $w$ can be directly derived from rotations on the circle as
\begin{equation}\label{eq:floor}
w_n=s_{\alpha, \rho}(n)=[(n+1)\alpha+\rho]-[n\alpha+\rho],\quad n=0,1,2,\dots.
\end{equation} or as
\begin{equation}\label{eq:ceil}
w_n=s'_{\alpha, \rho}(n)=\lceil(n+1)\alpha+\rho\rceil-\lceil n\alpha+\rho\rceil,\quad n=0,1,2,\dots.
\end{equation}
Here  $0<\alpha<1$  and $\rho$  are real numbers, $[\cdot]$ is the floor function, and $\lceil \cdot \rceil$ is the ceiling function.

Sturmian words have been named after Jacques Charles Fran\c{c}ois Sturm, who never studied them. A whole chapter is dedicated to them in Lothaire's book `Algebraic combinatorics on words' (\cite{Lothaire}).  There is a huge literature, in particular on the \emph{homogeneous} Sturmian words $$c_\alpha:=s_{\alpha, \alpha},$$    which have been studied since Johann III Bernoulli. The homogeneous  Sturmian words are also known as \emph{characteristic words} see Chapter 9 in \cite{Allou}.

Interestingly, for certain $\alpha$ and $\rho$  the Sturmian word $w$ is a fixed point $\sigma(w)=w$ of a morphism\footnote{We interchangeably use the terms  morphisms and  substitutions.} $\sigma$ of the monoid of words over the alphabet $\{0,1\}$. For example, for $\alpha=\rho=(3-\sqrt{5})/2$ one obtains the Fibonacci word $c_\alpha=0100101\dots$, fixed point of the Fibonacci morphism $\varphi$ given by $\varphi(0)=01$, $\varphi(1)=0$. Another example is the Pell word $c_{\alpha}=0010010001\dots$ obtained for $\alpha=\rho=(2-\sqrt{2})/2$, with morphism given by $0\rightarrow 001,\; 1\rightarrow 0$. 

It is well known for which $\alpha$ one obtains a morphism invariant $c_{\alpha}$. This was first obtained in \cite{Moran}, and an extensive treatment can be found in \cite[Section 2.3.6]{Lothaire}. The result is that $\alpha\in (0,\tfrac12)$ gives a fixed point if and only if there exists a natural number $k$ such that $\alpha$ has continued fraction expansion
\begin{equation}\label{eq:small} \alpha=[0;1+a_0,\overline{a_1\dots a_k}], \qquad a_k\ge a_0\ge 1, \end{equation}
and $\alpha\in (\tfrac12,1)$ gives a fixed point if and only if there exists a natural number $k$ such that $\alpha$ has continued fraction expansion
\begin{equation}\label{eq:large} \alpha=[0;1,a_0,\overline{a_1\dots a_k}], \qquad a_k\ge a_0.\end{equation}
The Fibonacci word is obtained for $k=1, a_0=a_1=1$, and the Pell word for $k=1, a_0=2, a_1=3$.

Any $\alpha$ that gives a substitution invariant $c_\alpha$ is called a \emph{Sturm number}. In terms of their continued fraction expansions these are characterized in equations (\ref{eq:small}) and (\ref{eq:large}). There is however a simple algebraic way to describe them, given in \cite{Allauzen}:

 an irrational  number $\alpha \in (0,1)$ is a Sturm number if and only if it is a quadratic irrational number whose algebraic conjugate $\overline{\alpha}$,  defined by the equation $(x-\alpha)(x-\overline{\alpha})=0$,  satisfies
 $$\overline{\alpha} \notin [0,1].$$
\indent A simple manipulation shows that for $\alpha\in (0,\tfrac12)$ the number $\beta =1-\alpha$ has an expansion as in equation (\ref{eq:large}) with the same $k$ and $a_j,\, j=0,\dots,k$. Moreover, the Sturmian word $w(\beta)$ is equal to the word $E(w(\alpha))$, where $E$ is the `exchange' morphism
$$E: \; {\mor{1}{0} }\:.$$
The latter is shown in the proof of Theorem 2.3.25 in \cite{Lothaire}. Note that this implies that if $\sigma$ generates $w(\alpha)$, then $E\sigma E$ generates $w(1-\alpha)$. Because of this duality we will confine ourselves often to $\alpha$ with $0<\alpha<\tfrac12$ in the sequel.

The first question we will consider is: what are the morphisms that leave a homogeneous Sturmian sequence $c_\alpha$ invariant? The answer in \cite{Moran} is: they are compositions of the infinitely many morphisms $G_k:\; 0\rightarrow 1^k0,\, 1\rightarrow 1$ and $H_k=G_kE$. The answer in \cite{Allou} is: they are compositions of the infinitely many morphisms $h_k:\; 0\rightarrow 0^k1,\, 1\rightarrow 0\rightarrow 0^k10$ (actually only for $\alpha$'s with a purely periodic continued fraction expansion). See \cite{Komatsu} for yet another infinite family of morphisms.

 In the paper \cite{ItoYasutomi} the authors call the inhomogeneous sequence $s_{\alpha, 0}$ a characteristic sequence, and do actually derive a result close to our Theorem \ref{th:diam0}, using completely different techniques with continued fractions and extensive matrix multiplications.

More satisfactory is the answer in the book \cite{Lothaire} or the paper \cite{Berstel}, where only two generating morphisms are used, namely the exchange morphism $E$ and the morphism $G$ given by $G(0)=:0,\, G(1)=01.$ What we propose are also only two generators, which we denote by $\varphi_0$ and $\varphi_1$, given by
$$\varphi_0: \; \mor{0}{01} \qquad   \varphi_1: \; \mor{01}{0}\:.$$
Note that $\varphi_0=G$, and that $\varphi_1=GE$, the Fibonacci morphism. Obviously, this proposal is very close to the one in \cite{Lothaire}, but what we gain is a natural way to index all the morphisms that leave homogeneous Sturmian words invariant by a binary tree (actually two binary trees, one for $\alpha\in (0,\tfrac12)$, and a dual version for $\alpha\in (\tfrac12,1)$). In Section~\ref{sec:Kep} we treat some preliminaries to give in Section~\ref{sec:mor} our main result.

We remark that a similar tree associated to the rational numbers appears in the work of de Luca \cite{Luca1,Luca}. The labeling there is not with morphisms, but with words.

The second question we will consider is: what are the substitution invariant Sturmian words that can only be obtained via the ceiling function, i.e., the  Sturmian words that can only be obtained as in equation (\ref{eq:ceil})? In this respect the homogeneous Sturmian words are regular, in that for all $\alpha$   $$c_\alpha=s_{\alpha,\alpha}=s'_{\alpha,\alpha}.$$
So these `strictly ceiling' Sturmian words have to be sought among the inhomogeneous Sturmian words, what we do in Section \ref{sec:inhom}.

The short Section \ref{sec:gen} is more or less independent of the remainder of the paper, but its contents have been very useful in our research.


\section{Homogeneous Sturmian words}

\subsection{The binary tree of harmonic fractions}\label{sec:Kep}
 The binary tree is a graph with $2^n$ nodes $i_1\dots i_n$ at level $n$ for $n=1,2,\dots$, where the $i_k$ are 0 or 1. At level 0 there is the root node $\Lambda$.

\begin{wrapfigure}{R}{0.45\textwidth}
\vspace*{-.2cm}
\begin{center}
\begin{tikzpicture}
[level distance=8mm,
every node/.style={fill=magenta!10,circle,inner sep=1pt},
level 1/.style={sibling distance=35mm,nodes={fill=magenta!10}},
level 2/.style={sibling distance=18mm,nodes={fill=magenta!10}},
level 3/.style={sibling distance=8mm,nodes={fill=magenta!10}}]
\node {$\frac{1}{2}$}            
child {node {$\frac{1}{3}$}      
 child {node {$\frac{1}{4}$}     
  child {node {$\frac{1}{5}$}}   
  child {node {$\frac{4}{5}$}}   
}
 child {node {$\frac{3}{4}$}     
  child {node {$\frac{3}{7}$}}   
  child {node {$\frac{4}{7}$}}   
}
}
child {node {$\frac{2}{3}$}      
 child {node {$\frac{2}{5}$}     
  child {node {$\frac{2}{7}$}}   
  child {node {$\frac{5}{7}$}}   
}
 child {node {$\frac{3}{5}$}     
   child {node {$\frac{3}{8}$}}               
   child {node {$\frac{5}{8}$}}              
}
};
\end{tikzpicture}
\end{center}
\end{wrapfigure}

 As early as 1619 Johannes Kepler defined in \cite{Kepler} a binary tree with fractions $\tfrac{p}{q}$ at the nodes. In the root there is $\tfrac12$, and if  $\tfrac{p}{q}$ is at a node, then the  two children nodes receive the fractions\\[-.3cm]

$\hspace*{1cm}    \dfrac{p}{p+q},\quad \dfrac{q}{p+q}.$\\

  Rather surprisingly, every rational number $p/q$ with $(p,q)=1$ in the interval (0,1) occurs exactly once in the tree. This is not hard to prove, see, e.g., the paper \cite{Smolinsky}. We remark that the paper \cite{Kimber} consider this problem for larger classes of  trees, but regretfully the rules for what the author's call the Kepler tree are different from Kepler's, but rather like those for the Calkin-Wilf tree (see \cite{Calkin}).

\noindent We introduce the two $2\times 2$ matrices\\[-.2cm]

$\hspace*{1cm} K_0:=\left( \begin{matrix} 1\, 0\\ 1\, 1 \end{matrix} \right), \quad K_1:=\left( \begin{matrix} 0\, 1\\ 1\, 1 \end{matrix} \right),$\\

\noindent which we call the \emph{Kepler matrices}.

\begin{wrapfigure}{R}{0.53\textwidth}
\vspace*{-2.0cm}
\begin{center}
\begin{tikzpicture}
[level distance=8mm,
every node/.style={fill=blue!10,rectangle,inner sep=1pt},
level 1/.style={sibling distance=40mm,nodes={fill=blue!10}},
level 2/.style={sibling distance=20mm,nodes={fill=blue!10}},
level 3/.style={sibling distance=10mm,nodes={fill=blue!10}}]
\node {\MM{1}{0}{0}{1}}            
child {node {\MM{1}{0}{1}{1}}      
 child {node {\MM{1}{0}{2}{1}}     
  child {node { \MM{1}{0}{3}{1}}}   
  child {node {\MM{2}{1}{3}{1}}}   
}
 child {node {\MM{1}{1}{2}{1}}     
  child {node {\MM{1}{1}{3}{2}}}   
  child {node {\MM{2}{1}{3}{2}}}   
}
}
child {node {\MM{0}{1}{1}{1}}      
 child {node {\MM{0}{1}{1}{2}}     
  child {node {\MM{0}{1}{1}{3}}}   
  child {node {\MM{1}{2}{1}{3}}}   
}
 child {node {\MM{1}{1}{1}{2}}     
   child {node {\MM{1}{1}{2}{3}}}  
   child {node {\MM{1}{2}{2}{3}}}  
}
};
\end{tikzpicture}
\end{center}
\end{wrapfigure}

It is clear that the fraction at the node $\underline{i}=i_1\dots i_n$ in Kepler's tree of fractions is equal to $p/q$,
where\\[-.2cm]

$\hspace*{1cm}\PQ{p\,}{q\, }= K_{i_n}\cdots K_{i_1}\PQ{1}{2}.$\\

We claim that all the matrices $K_{i_n}\cdots K_{i_1}$ are different when $n$ ranges over the natural numbers, and $i_1\dots i_n$ is a string of 0's and 1's. Formulated slightly differently we have the following.

\begin{lemma} \label{lem:free} The monoid of matrices generated by $K_0$ and $K_1$ is free.
\end{lemma}

\noindent \emph{Proof:} If $K_{\underline{i}}=K_{\underline{j}}$, then $K_{\underline{i}}\pq{1}{2}=K_{\underline{j}}\pq{1}{2}$, contradicting  uniqueness on the Kepler tree. \fbox

\medskip

We remark that in general it is hard to determine freeness of matrix monoids. It is for instance an undecidable problem for $3\times 3$ nonnegative integer matrices (\cite{Klarner}, see also \cite{Cassaigne}). We mention also that $K_0$ and $K_1$ are unimodular matrices, but that they do not satisfy the criteria in \cite{Lyndon}.


\subsection{A tree of  morphisms}\label{sec:mor}

Let  $\varphi_0$ and $\varphi_1$ be the two morphisms given  by
$$\varphi_0: \; \mor{0}{01} \qquad   \varphi_1: \; \mor{01}{0}.$$
We form a tree of morphisms $\mathcal{T}_\varphi$ by putting ${\rm Id}\!:0\rightarrow 0,\, 1\rightarrow 1$ at $\Lambda$, and $\varphi_{i_n}\cdots \varphi_{i_1}$ at node $\underline{i}=i_1\dots i_n$ for all $n$ and all $i_k \in \{0,1\}$, $k=1,\dots,n$.

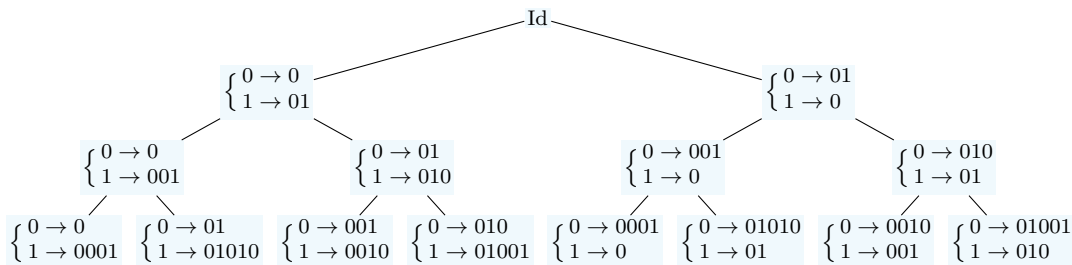
\smallskip

\begin{center}
\begin{tikzpicture}
[level distance=10mm,
every node/.style={fill=cyan!5,rectangle,inner sep=1pt},
level 1/.style={sibling distance=72mm,nodes={fill=cyan!5}},
level 2/.style={sibling distance=36mm,nodes={fill=cyan!5}},
level 3/.style={sibling distance=18mm,nodes={fill=cyan!5}}]
\node {\footnotesize  Id}            
child {node {\footnotesize $\mor{0}{01}$}            
 child {node {\footnotesize $\mor{0}{001}$}          
  child {node {\footnotesize $\mor{0}{0001}$}}       
  child {node {\footnotesize $\mor{01}{01010}$}}     
}
 child {node {\footnotesize $\mor{01}{010}$}         
  child {node {\footnotesize $\mor{001}{0010}$}}     
  child {node {\footnotesize $\mor{010}{01001}$}}    
}
}
child {node {\footnotesize $\mor{01}{0}$}            
 child {node {\footnotesize $\mor{001}{0}$}          
  child {node {\footnotesize $\mor{0001}{0}$}}       
  child {node {\footnotesize $\mor{01010}{01}$}}     
}
 child {node {\footnotesize $\mor{010}{01}$}         
   child {node {\footnotesize $\mor{0010}{001}$}}    
   child {node {\footnotesize $\mor{01001}{010}$}}   
}
};
\end{tikzpicture}
\end{center}

\medskip

The figure shows the first 3 levels of this tree, labeled with the morphisms. Note that the left edge of $\mathcal{T}_\varphi$ with nodes $\underline{i}=0^n$ contains the morphisms $\varphi_0^n$, which do not generate infinite words.

\begin{theorem} \label{th: hom} The tree $\mathcal{T}_\varphi$  contains all morphisms that have homogeneous Sturmian words $c_{\alpha}$ as fixed point, for any $\alpha$  with $0<\alpha<\tfrac12$. Each such morphism occurs exactly once.
\end{theorem}

\emph{Proof:} In \cite{Lothaire} it is proved that for $\alpha\in (0,1)$, any morphism $f$ fixing  a homogeneous Sturmian word is a composition of the two morphisms $E$ and $G$, excluding $f=E^n$, $f=G^n$ and $f=EG^nE$ for $n\ge 1$. Moreover, if $0<\alpha<\tfrac12$, then the first element in the composition of $f$ is $G$. But since $E^2={\rm Id}$, $f$ can then be written as a composition of $G=\varphi_0$ and $GE=\varphi_1$. This finishes the existence part of the proof.

For the uniqueness part, we remark first that it is shown in \cite[Corollary 2.3.15]{Lothaire},  that  in the  monoid generated by the two morphisms $E$ and $GE$ the only relation  is $E^2={\rm Id}$ . This implies, of course, that the monoid generated by $G=\varphi_0$ and $GE=\varphi_1$, is free, but here we prefer to give a short self-contained proof, proving something stronger, which yields the emergence of the binary tree.

Consider the incidence  matrices of the morphisms $\varphi_0$ and $\varphi_1$ :\\[-.2cm]

$\hspace*{1cm}  M_0:=\left( \begin{matrix} 1\, 1\\ 0\, 1 \end{matrix} \right), \quad M_1:=\left( \begin{matrix} 1\, 1\\ 1\, 0 \end{matrix} \right).$\\

\begin{wrapfigure}{R}{0.45\textwidth}
\vspace*{-0.8cm}
\begin{center}
\begin{tikzpicture}
[level distance=8mm,
every node/.style={fill=violet!5,rectangle,inner sep=1pt},
level 1/.style={sibling distance=35mm,nodes={fill=violet!5}},
level 2/.style={sibling distance=18mm,nodes={fill=violet!5}},
level 3/.style={sibling distance=8mm,nodes={fill=violet!5}}]
\node {\MM{1}{0}{0}{1}}            
child {node {\MM{1}{1}{0}{1}}      
 child {node {\MM{1}{2}{0}{1}}     
  child {node {\MM{1}{3}{0}{1}}}   
  child {node {\MM{1}{3}{1}{2}}}   
}
 child {node {\MM{1}{2}{1}{1}}     
  child {node {\MM{2}{3}{1}{1}}}   
  child {node {\MM{2}{3}{1}{2}}}   
}
}
child {node {\MM{1}{1}{1}{0}}      
 child {node {\MM{2}{1}{1}{0}}     
  child {node {\MM{3}{1}{1}{0}}}   
  child {node {\MM{3}{1}{2}{1}}}   
}
 child {node {\MM{2}{1}{1}{1}}     
   child {node {\MM{3}{2}{1}{1}}}  
   child {node {\MM{3}{2}{2}{1}}}  
}
};
\end{tikzpicture}
\end{center}
\end{wrapfigure}

Obviously morphisms with different incidence matrices are different. Let ${\mathcal{T}}_M$ be the tree with the matrix product $M_{i_n}\cdots M_{i_1}$ at node $i_1\dots i_n$.

The matrices $M_0$ and $M_1$ are conjugate to the matrices $K_0$ and $K_1$ by the same conjugation matrix \MM{0}{1}{1}{0}.
It follows that for \emph{any}  node $\underline{i}$ one has \\[-.2cm]

$\hspace*{1cm} M_{\underline{i}}= \left( \begin{matrix} 0\, 1\\ 1\, 0 \end{matrix} \right) K_{\underline{i}} \left( \begin{matrix} 0\, 1\\ 1\, 0 \end{matrix} \right).$\\

\noindent But then Lemma~\ref{lem:free} implies that all the $M_{\underline{i}}$ on ${\mathcal{T}}_M$ are different, and so each morphism occurs exactly once on $\mathcal{T}_\varphi$. \fbox

\subsection{A tree of Sturm numbers}

The tree of morphisms $\mathcal{T}_\varphi$ has a left edge with morphisms that do not generate infinite words.
  Below we display the  first three levels of the tree of Sturm numbers $\alpha$ with $\alpha<\frac12$ associated to the morphisms of  $\mathcal{T}_\varphi$. Each such $\alpha$ will occur infinitely times, since the powers of a morphism generate the same Sturmian word. In particular we will see on the right edge the number $(3-\sqrt{5})/2$ associated to the the powers of the Fibonacci morphism  $\varphi_1^n$.



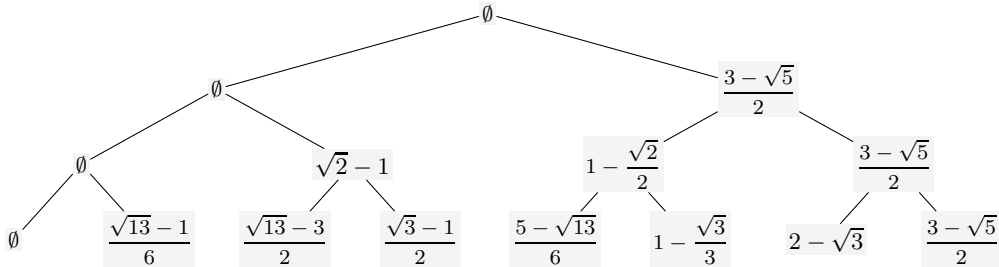
\begin{center}
\begin{tikzpicture}
[level distance=10mm,
every node/.style={fill=grey!5,rectangle,inner sep=1pt},
level 1/.style={sibling distance=72mm,nodes={fill=grey!5}},
level 2/.style={sibling distance=36mm,nodes={fill=grey!5}},
level 3/.style={sibling distance=18mm,nodes={fill=grey!5}}]
\node {\footnotesize $\emptyset$}            
child {node {\footnotesize $\emptyset$}            
 child {node {\footnotesize $\emptyset$}          
  child {node {\footnotesize $\emptyset$}}       
  child {node {\footnotesize $\dfrac{\sqrt{13}-1}{6}$}}     
}
 child {node {\small $\sqrt{2}-1$}         
  child {node {\footnotesize $\dfrac{\sqrt{13}-3}{2}$}}     
  child {node {\footnotesize $\dfrac{\sqrt{3}-1}{2}$}}    
}
}
child {node {\footnotesize $\dfrac{3-\sqrt{5}}{2}$}            
 child {node {\footnotesize $1-\dfrac{\sqrt{2}}{2}$}          
  child {node {\footnotesize $\dfrac{5-\sqrt{13}}{6}$}}       
  child {node {\footnotesize $1-\dfrac{\sqrt{3}}{3}$}}     
}
 child {node {\footnotesize $\dfrac{3-\sqrt{5}}{2}$}         
   child {node {\small $2-\sqrt{3}$}}    
   child {node {\footnotesize $\dfrac{3-\sqrt{5}}{2}$}}   
}
};
\end{tikzpicture}
\end{center}


\section{Inhomogeneous Sturmian words}\label{sec:inhom}

  In this section we consider \emph{all} substitution invariant Sturmian words. There is again a simple algebraic characterization given by Yasutomi in  (\cite{Yasutomi35}): 

  Let $0 < \alpha < 1$ and $0 \le \rho \le 1$. Then $s_{\alpha,\rho}$ is
substitution invariant if and only if the following two conditions are satisfied:
\begin{eqnarray*}
&(i)& \alpha \;{\rm is\, an\, irrational\, quadratic\, number\, and\;} \rho\in \mathbb{Q}(\alpha);\\
&(ii)& \overline{\alpha} > 1,\, 1 - \overline{\alpha} \le \overline{\rho} \le\overline{\alpha}\quad{\rm or\quad} \overline{\alpha} < 0,\, \overline{\alpha} \le \overline{\rho} \le 1 - \overline{\alpha}.
\end{eqnarray*}

\subsection{The eight elementary  morphisms}

We define the morphisms $\psi_i$ for $i=1,\dots,8$ by
\begin{eqnarray*}  &\psi_1:& \; 0\rightarrow 01,\, 1\rightarrow 0, \quad \psi_2: \; 0\rightarrow 10,\, 1\rightarrow 0, \quad
\psi_3: \;\; 0\rightarrow 0,\, 1\rightarrow 01, \quad \psi_4: \;\; 0\rightarrow 0,\, 1\rightarrow 10.\\
 &\psi_5:& \; 0\rightarrow 1,\, 1\rightarrow 10, \quad \psi_6: \; 0\rightarrow 1,\, 1\rightarrow 01,\quad
\psi_7: \;\;0\rightarrow 10,\, 1\rightarrow 1 , \quad \psi_8: \;\; 0\rightarrow 01,\, 1\rightarrow 1.
\end{eqnarray*}
 In the previous section the two morphisms  $\psi_1=\varphi_1$, and $\psi_3=\varphi_0$ were used. The first four morphisms are linked to Sturmian words with slope $\alpha<1/2$, the last four to Sturmian words with slope $\alpha>1/2$. In the columns there is duality: $\psi_{i+4}=E\psi_i E$ for $i=1,2,3,4$.
 We also have $\psi_{2i}=\widetilde{\psi}_{2i-1}$ for $i=1,2,3,4$, where $\widetilde{\sigma}$ is the time reversal of a morphism $\sigma$.

 The notation used in \cite{Lothaire} is:\\[-.35cm]
$$\psi_1=\varphi,   \quad \psi_2= \widetilde{\varphi}, \quad \psi_3=\varphi E,   \quad \psi_4= \widetilde{\varphi}E,\quad
\psi_5=E\varphi E,   \quad \psi_6= E\widetilde{\varphi}E, \quad \psi_7=E\varphi,   \quad \psi_8= E\widetilde{\varphi}.$$\\[-.5cm]
Let $ \mathcal{M}_{i,j}=\langle \psi_i,\psi_j \rangle$ denote the monoid generated by the morphisms $\psi_i$ and $\psi_j$. We will also need $ \mathcal{M}_i=\langle \psi_i\rangle$, the set of powers of $\psi_i$.


\subsection{The floor-ceiling structure of Sturmian sequences}
 For most $\alpha$'s and $\rho$'s  the floor and the ceiling representation of a Sturmian word in equations (\ref{eq:floor}) and (\ref{eq:ceil}) are equal. Rather surprisingly, \emph{if} they are not equal, then they only differ in at most two consecutive indices (\cite{Lothaire}).
 If there exists a natural number $m_{\diamond}$ such that
 $$ s_{\alpha,\rho}(m_{\diamond}-1)\ne s'_{\alpha,\rho}(m_{\diamond}-1)\; {\rm and\;}  s_{\alpha,\rho}(m_{\diamond})\ne s'_{\alpha,\rho}(m_{\diamond}),$$
 then we call $(s_{\alpha,\rho}, s'_{\alpha,\rho})$ a \emph{lozenge pair} with index $m_{\diamond}$. In case $m_{\diamond}=0$, there is actually only the index 0 where they differ. As indicated on page 48 of \cite{Lothaire}, $(s_{\alpha,\rho}, s'_{\alpha,\rho})$ is a lozenge pair with index $m_{\diamond}$ if and only if
 \begin{equation}\label{eq:loz}  \alpha m_\diamond +\rho \in \mathbb{N}.\end{equation}

 \medskip

\noindent {\bf Example} Let $\alpha=(3-\sqrt{5})/2$ and $\rho=(\sqrt{5}-1)/2$. Then $\alpha+\rho=1$, so $(s_{\alpha,\rho}, s'_{\alpha,\rho})$ is a lozenge pair with index 1. Here $s_{\alpha,\rho}=1001001010010010\dots $ and $s'_{\alpha,\rho}=0101001010010010,\dots$. Both words are substitution invariant for the substitution $0\rightarrow 010, 1\rightarrow 10$.

\medskip

The following result is related to Corollary 1.4. in \cite{Berthe-et-al}.

\begin{proposition} \label{prop:loz} For substitution invariant lozenge pairs  $m_{\diamond}=0$ or $m_{\diamond}=1$.
\end{proposition}

\noindent  \emph{Proof:} This follows directly from equation (\ref{eq:loz}) and Yasutomi's characterization.\\
Suppose $m_\diamond \ge 2$ and $\alpha m_\diamond +\rho=k$ for an integer $k$. Since $0<\alpha m_{\diamond}<m_\diamond $ and $0\le \rho \le 1$ we must have
$$\alpha m_\diamond +\rho=k, {\rm where \;} k \in \{1,2,\dots,  m_\diamond \}.$$
One easily checks that  $\overline{\rho}=k-\overline{\alpha} m_\diamond$. Now if $\overline{\alpha}>1$, then it should hold that $$1-\overline
{\alpha}\le \overline{\rho}=k-  \overline{\alpha} m_\diamond\; \Rightarrow \; (m_\diamond-1)\overline{\alpha}\le k-1
\;\Rightarrow\; \overline{\alpha}\le 1,$$
yielding  a contradiction. Similarly the case $\overline{\alpha}<0$ yields a contradiction. \qed

\medskip

We undertake the task of  determining all substitution invariant Sturmian words that are lozenge pairs.
We start with  the case $m_{\diamond}=1$, which is simpler than $m_{\diamond}=0$.

\subsection{Substitution invariant Sturmian words with {$\mathbf{m_{\mathbf{\diamond}}=1}$}}

\quad
Note  that $m_\diamond=1$ implies
$$\alpha+\rho=1.$$
We first give a very simple way to obtain the lozenge pair.

\begin{proposition}\label{prop:loz01} For $\alpha\in (0,1)$  let  $(s_{\alpha, 1-\alpha},s'_{\alpha, 1-\alpha})$ be the lozenge pair with  $m_{\diamond}=1$. Then
$$s_{\alpha, 1-\alpha}=10\,c_\alpha\quad{\rm and\quad} s'_{\alpha, 1-\alpha}=01\,c_\alpha.$$
\end{proposition}

\noindent \emph{Proof:} Since $\alpha<1,\, s_{\alpha, 1-\alpha}(0)=[\alpha+\rho]-[\rho]=1$ and $s_{\alpha, 1-\alpha}(1)= [2\alpha+\rho]-[\alpha+\rho]=0$. Also, $s'_{\alpha, 1-\alpha}(0)=\lceil\alpha+\rho\rceil-\lceil\rho\rceil=0$ and $s'_{\alpha, 1-\alpha}(1)= \lceil2\alpha+\rho\rceil-\lceil\alpha+\rho\rceil=1$. Let $S$ be the shift: $S(w_0w_1w_2\dots)=w_1w_2\dots$. Adding $\alpha$ to $\rho$ shifts a Sturmian sequence by one: $s_{\alpha,\rho+\alpha}=S(s_{\alpha,\rho})$.  So $S^2(s_{\alpha, 1-\alpha})=s_{\alpha, 1+\alpha}=s_{\alpha, \alpha}=c_\alpha$. \qed

\medskip

We still have to investigate whether $s_{\alpha, 1-\alpha}$ and $s'_{\alpha, 1-\alpha}$ are substitution invariant. This can be directly derived from \cite{Berstel-Seebold}, but we give here a short and more global proof. We start with a combinatorial lemma.

\begin{lemma}\label{lem:pal}  For any $\psi\in  \mathcal{M}_{1,3}$ the words $01\psi^2(0)$ and $10\psi^2(1)$ are palindromes.
\end{lemma}

\noindent \emph{Proof:}  This relies on the notions and results of  \cite[Section 2.2.1]{Lothaire}. The standard morphisms are the elements of
$\langle \varphi,E\rangle=\langle \psi_{1},E \rangle$. Since $\psi_3=\varphi E$, all $\psi$ from $ \mathcal{M}_{1,3}$ are standard.
By Proposition 2.2.2 and Proposition 2.3.11 of \cite{Lothaire}, the words $\psi(0)$ and $\psi(1)$ are two standard words, which differ in their last two letters.
 But then $\psi^2(0)$ and $\psi^2(1)$ are  standard words so that $\psi^2(0)$ ends in 0, and $\psi^2(1)$ ends in 1. Moreover, according to \cite[Theorem 2.2.4]{Lothaire}, a word $w$ is standard if and only if it has length 1 or there exists a palindrome word $p$ such that $w=p\,01$ or $w=p\,10$. So the words $$01\psi^2(0)=01p\,10\quad {\rm and\;} 10\psi^2(1)=10 p'01$$
are palindromes. The length 1 case may occur, but then 010 is a palindrome. \qed

\bigskip

\begin{theorem} \label{th:diam1} Let $s_{\alpha, \rho}$ and $s'_{\alpha, \rho}$ be substitution invariant Sturmian words with $m_{\diamond}=1$ and $\alpha<1/2$.
Then these two words are fixed points of $\psi^2$, where $\psi \in \mathcal{M}_{2,4}\smallsetminus \mathcal{M}_{4}$.
Here $\psi$ is given by $\widetilde{\psi}(c_\alpha)=c_\alpha$.
\end{theorem}

\noindent \emph{Proof:} The condition $m_{\diamond}=1$ means that $\rho=1-\alpha$. What we will show is that $s_{\alpha, 1-\alpha}$ and $s'_{\alpha, 1-\alpha}$ are fixed points of $\widetilde{\psi}^2$, where $\psi$ fixes $c_\alpha$. Recall that $\widetilde{\psi}$ is the time reversal of $\psi$,
and that $\widetilde{\psi}_1=\psi_2$, and $\widetilde{\psi}_3=\psi_4$.
In general one has $\widetilde{\sigma \tau} = \widetilde{\sigma}\widetilde{\tau}$ for two morphisms $\sigma$ and $\tau$.
It thus follows from Lemma \ref{lem:pal} that for $\psi \in \mathcal{M}_{1,3}$ and all $n\ge 1$
$$01\,\psi^{2n}(0)=\widetilde{\psi}^{2n}(0)\,10\quad {\rm and\quad} 10\,\psi^{2n}(1)=\widetilde{\psi}^{2n}(1)\,01.$$
For such $\psi$, not from $\mathcal{M}_{3}$, when $n\rightarrow\infty$,  the left sides converge to $01c_{\alpha}$, respectively $10c_{\alpha}$, and thus by Proposition \ref{prop:loz01}, the right sides converge to $s_{\alpha, 1-\alpha}$ respectively $s'_{\alpha, 1-\alpha}$. \qed

\medskip

 It is easily seen via duality that Theorem \ref{th:diam1} also applies in case $\alpha>1/2$, where $\mathcal{M}_{2,4}\smallsetminus \mathcal{M}_{4}$ has to be replaced by $\mathcal{M}_{6,8}\smallsetminus \mathcal{M}_{8}$.

\subsection{Substitution invariant Sturmian words with {$\mathbf{m_{\mathbf{\diamond}}=0}$} }

\quad
Note  that $m_\diamond=0$ implies $\rho=0.$
As in the case of $m_\diamond=1$ there is a simple way to obtain the lozenge pair:
$$s_{\alpha, 0}=0\,c_\alpha\quad{\rm and\quad} s'_{\alpha, 0}=1\,c_\alpha.$$
It is well known that $0\,c_\alpha$ is substitution invariant, see, e.g., Corollary 1.4. in \cite{Berthe-et-al}.
However, it is  now less simple to determine the substitution fixing $0\,c_\alpha$ for a given $\alpha$.

\medskip

\noindent {\bf Example} Let $\alpha=(\sqrt{13}-1)/6$. Then $\gamma(c_\alpha)=c_\alpha$ for $\gamma$ given by
$\gamma(0)=01,\; \gamma(1)=01010.$\\(The same morphism is considered in \cite{Berstel-Seebold} on page 262.)\\
Here $\psi(s_{\alpha, 0})=s_{\alpha, 0}$ for $\psi$  given by\\[0.1cm]
\hspace*{4cm}$\psi(0)=0010101,\quad \psi(1)=0010101001010101.$

\medskip


A recipe is given in \cite{Berstel-Seebold}. The recipe depends strongly on the last letter of $\gamma(0)$, so it is useful to characterize the morphisms $\gamma$ with $\gamma(0)$ ending in $0$.

\medskip

\begin{lemma} Let $\gamma=\psi_{i_1}\dots\psi_{i_m}$ be a morphism from $\mathcal{M}_{1,3}$, then $\gamma(0)$ ends in $0$ if
and only if the number of $1$ in $i_1\dots i_m$ is even.
\end{lemma}

\smallskip

\noindent {\emph{Proof:}} All three words $\psi_1(1), \psi_3(0),\psi_3(1)$ end in $0$, but $\psi_1(0)$ ends in $1$. \hfill\qed

\medskip

\noindent We denote the set of $\gamma$ from $\mathcal{M}_{1,3}$ such  $\gamma(0)$ ends in $0$ by $\mathcal{M}^0_{1,3}$.

\begin{proposition}\label{prop:psi-g}
 Let $\gamma \in \mathcal{M}^0_{1,3}$, such that $\gamma(c_\alpha)=c_\alpha$.  Let $\Psi_\gamma$ be conjugate to $\gamma$, with conjugating word
 equal to $u=\gamma(0)0^{-1}$. Then $$\Psi_{\gamma}(0c_\alpha) = 0c_\alpha.$$
\end{proposition}

\noindent For a  proof of this proposition, see \cite[Theorem 3.1]{Berstel-Seebold}.

\medskip

\begin{proposition}\label{prop:psi-g-d}
 Let $\gamma,\delta \in \mathcal{M}^0_{1,3}$.  Then $\Psi_{\gamma\delta} = \Psi_{\gamma}\Psi_{\delta}.$
\end{proposition}

\noindent {\emph{Proof:}} We assume that $|\gamma(0)|<|\gamma(1)|$, the proof of the other case is quite similar.
It is well known that there exist words $u,v,x,$ and $y$ such that
\begin{equation*}\label{eq:gd}
\gamma : \left\{ \begin{aligned}
0 & \rightarrow u\,0\\
1 & \rightarrow u\,0\,v
\end{aligned}\right.
\quad \Psi_\gamma : \left\{ \begin{aligned}
0 & \rightarrow 0\,u\\
1 & \rightarrow 0\,v\,u
\end{aligned}\right.
\qquad \delta : \left\{ \begin{aligned}
0 & \rightarrow x\,0,\\
1 & \rightarrow x\,0\,y
\end{aligned}\right.
\quad \Psi_\delta : \left\{ \begin{aligned}
0 & \rightarrow 0\,x\\
1 & \rightarrow 0\,y\,x.
\end{aligned}\right.
\end{equation*}
The product $\gamma\delta$ also generates a characteristic sequence, and we have
\begin{equation*}\label{eq:gd}
\gamma\delta : \left\{ \begin{aligned}
0 & \rightarrow \gamma(x)\,\gamma(0)=\gamma(x)\,u\,0\\
1 & \rightarrow \gamma(x)\,\gamma(0)\,\gamma(y)=\gamma(x)\,u\,0\,\gamma(y)
\end{aligned}\right.
\qquad \Psi_{\gamma\delta} : \left\{ \begin{aligned}
0 & \rightarrow 0\,\gamma(x)\,u\\
1 & \rightarrow 0\,\gamma(y\,)\gamma(x)\,u.
\end{aligned}\right.
\end{equation*}
This is a slightly extended version of \cite[Lemma 2.3.17 (iii)]{Lothaire}, which leads to
$\Psi_{\gamma\delta} = \Psi_{\gamma}\Psi_{\delta},$
since\begin{equation*}\label{eq:gd}
\Psi_{\gamma}\Psi_{\delta} : \left\{ \begin{aligned}
0 & \rightarrow \Psi_\gamma(0)\,\Psi_\gamma(x)=0\,u\,\Psi_\gamma(x)=0\,\gamma(x)\,u\\
1 & \rightarrow \Psi_\gamma(0)\,\Psi_\gamma(y)\,\Psi_\gamma(x)=0\,u\,\Psi_\gamma(y)\,\Psi_\gamma(x)=0\,\gamma(y\,)\gamma(x)\,u,
\end{aligned}\right.
\end{equation*}
where we use the conjugation relation $u\,\Psi_\gamma(w)=\gamma(w)\,u$ for $w=x, yx$.\hfill\qed

\medskip

\noindent For a further description we need next to $\psi_3$ yet another elementary morphism $\psi_8=E\,\widetilde{\varphi}$:
$$\psi_3: \; 0\rightarrow 0,\, 1\rightarrow 01,\qquad\psi_8: \; 0\rightarrow 01,\, 1\rightarrow 1.$$

\medskip

\begin{lemma}\label{lem:131} Let $\gamma=\psi_{1}\psi_{3}^n\psi_{1}$ for some $n\ge 0$. Then  $\Psi_\gamma=\psi_3\psi_8^{n+1}$.
\end{lemma}

\noindent {\emph{Proof:}} One easily finds that for all $n\ge0$
 $$\psi_{1}\psi_{3}^n\psi_{1}(0)=(01)^{n+1}0,\;\psi_{1}\psi_{3}^n\psi_{1}(1)=01;\quad \psi_3\psi_8^{n+1}(0)=0(01)^{n+1},\; \psi_3\psi_8^{n+1}(1)=01.$$
 This implies the statement of the lemma.\hfill\qed

 \medskip

\noindent We need more details on the structure of the map $\gamma\mapsto\Psi_\gamma$.

 \begin{proposition}\label{prop:p3p5} Let $\gamma$ be a morphism from $\mathcal{M}^0_{1,3}$. Then
 $$\Psi_{E\gamma E}=E\,\widetilde{\Psi}_\gamma\,E.$$
 Moreover, $\Psi_{E\gamma E}=\Psi_\gamma^*$, where $\cdot^*$ is the homomorphism defined by $\psi_3^*=\psi_8, \, \psi_8^*=\psi_3$.
\end{proposition}

\noindent {\emph{Proof:}} We assume that $|\gamma(0)|<|\gamma(1)|$, the proof of the other case is quite similar. We use the same facts on $\gamma$ from \cite{Lothaire} as in the proof of Lemma \ref{lem:pal}. Since $(\gamma(0), \gamma(1))$ is a standard pair, and all $\gamma(0)$ from $\mathcal{M}_{1,3}$ start with 0, there exist\footnote{Here we admit $p=0^{-1}$, with length $|p|=-1$.} palindromes $p$ and $q$ such that
\begin{equation*}
\gamma : \left\{ \begin{aligned}
0 & \rightarrow 0\,p\,0\,1\,0\\
1 & \rightarrow 0\,p\,0\,1\,0\,q\,0\,1,
\end{aligned}\right.
\qquad \Psi_{\gamma} : \left\{ \begin{aligned}
0 & \rightarrow 0\,0\,p\,0\,1\\
1 & \rightarrow 0\,q\,0\,1\,0\,p\,0\,1.
\end{aligned}\right.
\end{equation*}
So we have
\begin{equation*}
\widetilde{\Psi}_\gamma : \left\{ \begin{aligned}
0 & \rightarrow 1\,0\,p\,0\,0\\
1 & \rightarrow 1\,0\,p\,0\,1\,0\,q\,0,
\end{aligned}\right.
\qquad E\widetilde{\Psi}_\gamma E : \left\{ \begin{aligned}
0 & \rightarrow 0\,1\,\overline{p}\,1\,0\,1\overline{q}\,1\\
1 & \rightarrow 0\,1\,\overline{p}\,1\,1.
\end{aligned}\right.
\end{equation*}
On the other hand,
\begin{equation*}
E\gamma E : \left\{ \begin{aligned}
0 & \rightarrow 1\,\overline{p}\,1\,0\,1\overline{q}\,1\,0\\
1 & \rightarrow  1\,\overline{p}\,1\,0\,1,
\end{aligned}\right.
\qquad \Psi_{E\gamma E} : \left\{ \begin{aligned}
0 & \rightarrow  0\,1\,\overline{p}\,1\,0\,1\overline{q}\,1\\
1 & \rightarrow \dots.
\end{aligned}\right.
\end{equation*}
Note that we showed $\Psi_{E\gamma E}(0)=E\widetilde{\Psi}_\gamma E(0)$, but for 1 this is trickier.
A way to look at conjugation is by simultaneous repeated rotation. Here by \emph{rotation} we mean the map $\rho$ on words defined by
$$\rho(w_1w_2\dots w_m)=w_2\dots w_m w_1.$$
Moreover, words may only be rotated simultaneously if their first letters are equal.

\noindent From this viewpoint, $$\Psi_{E\gamma E}(0)=\rho^{|E\gamma E(0)|-1}E\gamma E(0)=\rho^{|p|+|q|+5}E\gamma E(0).$$
Since the length of $E\gamma E(1)$ equals $|p|+4$, and $E\gamma E(1)$ is a prefix of $E\gamma E(0)$, we may first rotate $|p|+4$ times, obtaining
$$\rho^{|p|+4}\Psi_{E\gamma E}(0)= \overline{q}\,1\,0\,1\,\overline{p}\,1\,0\,1  \quad \rho^{|p|+4}\Psi_{E\gamma E}(1)=1\,\overline{p}\,1\,0\,1.$$
To continue rotating $|q|+1$ times, , we must see that $\overline{q}\,1$ is a prefix of $1\,\overline{p}\,1\,0\,1$, \emph{and} we want to see that  the outcome is
$$\rho^{|q|+1}(1\,\overline{p}\,1\,0\,1)= 0\,1\,\overline{p}\,1\,1,\quad{\rm or\, equivalently,\quad}  \rho^{|q|+1}(0\,p\,0\,1\,0)= 1\,0\,p\,0\,0.$$
This requires for the time being that $|q|\le |p|+4$.

According to \cite[Theorem 2.2.4]{Lothaire} the word $0\,p\,0\,1\,0\,q$ is a palindrome, and so
$$0\,p\,0\,1\,0\,q=q\,0\,1\,0\,p\,0,$$
hence $\rho^{|q|+1}(0\,p\,0\,1\,0)$ has prefix $1\,0\,p\,0$, and the next letter is a $0$, because we see also from this equation that the first letter of $q$ is a $0$. Note that this palindrome equation also  implies that $\overline{q}\,1$ is a prefix of $1\,\overline{p}\,1\,0\,1$.

In the case that $|q|> |p|+4$, suppose that $|q|=a(|p|+4)+b$, where $b<|p|+4$. Then one rotates with $\rho^{|p|+4}$ $a+1$ times, followed by a rotation $\rho^{b+1}$. That this is possible,  one can deduce from  \cite[Proposition 2.2.2]{Lothaire}, which states that only the last two letters of the words $\gamma(0)\gamma(1)$ and $\gamma(1)\gamma(0)$ are different. A zig-zag argument then gives that $\gamma(1)=\gamma(0)^{a+1}01$.

For the second part of the proposition, note that $\psi_3^*=E\widetilde{\psi}_3E=\psi_8$, and $\psi_8^*=E\widetilde{\psi}_8E=\psi_3$.
Let $\Psi_\gamma=\psi_{i_1}\dots\psi_{i_m}$. Then we have, using the first part of the proposition,
$$\Psi_{E\gamma E}=E\,\widetilde{\Psi}_\gamma\,E=E\widetilde{\psi}_{i_1}\dots\widetilde{\psi}_{i_m}E=
E\widetilde{\psi}_{i_1}EE\widetilde{\psi}_{i_2}E\dots E\widetilde{\psi}_{i_m}E=\psi_{i_1}^*\psi_{i_2}^*\dots\psi_{i_m}^*=\Psi_\gamma^*. \hfill \qed$$

\begin{theorem}\label{th:diam0} Let $\alpha$ be a Sturmian number, with $0<\alpha<1$.  Then $s_{\alpha, 0}$ is a fixed point of some
 $\psi \in \mathcal{M}_{3,8}$.
 Conversely, any  $\psi \in \mathcal{M}_{3,8}\smallsetminus \{\mathcal{M}_3 \cup  \mathcal{M}_8\}$ fixes an  $s_{\alpha, 0}$. The same statements hold for $s'_{\alpha, 0}$, but then with $\mathcal{M}_{3,8}$ replaced by $\mathcal{M}_{4,7}$.
\end{theorem}

\noindent \emph{Proof:} We have $s_{\alpha, 0}=0\,c_\alpha$. Suppose $\gamma\in \mathcal{M}_{1,3}$ satisfies $\gamma(c_\alpha)=c_\alpha$. Then $\gamma^2(c_\alpha)=c_\alpha$ and $\gamma^2 \in \mathcal{M}^0_{1,3}$, so by Proposition \ref{prop:psi-g}, $s_{\alpha, 0}$ is fixed point of  $\Psi_{\gamma^2}$.

  We claim that any  $\Psi_{\gamma}$, where $\gamma$ is from $\mathcal{M}^0_{1,3}$,
 is an element of $\mathcal{M}_{3,8}$.
  We prove this claim by induction on  $m$ where $\gamma=\psi_{i_1}\dots\psi_{i_{m}}$.
  For $m= 2$,  $\gamma=\psi_1^2$, and the claim is true by Lemma \ref{lem:131}. Suppose the claim is true for all $\gamma$ from $\mathcal{M}^0_{1,3}$
  with length $m$ or less.
  An arbitrary $\gamma=\psi_{i_1}\dots\psi_{i_{m+1}}$ from $\mathcal{M}^0_{1,3}$, can be written as  $\gamma=\gamma'\gamma''$,
  where $\gamma'$ and $\gamma''$ are non-trivial elements of $\mathcal{M}^0_{1,3}$, \emph{unless} $\gamma$ has the form
  $$\gamma=\psi_{1}\psi_3^{m-1}\psi_{1},$$
  but then $\Psi_\gamma\in \mathcal{M}_{3,8}$ according to Lemma \ref{lem:131}.
  The first part of the theorem is proved.

  For the second part, we divide  the morphisms in $ \mathcal{M}_{3,8}$ into two types: the ones starting with $\psi_3$ and the ones starting with $\psi_8$. A density argument shows that first type corresponds to $\psi$ with $\alpha< 1/2$, and the second type to $\psi$ with $\alpha> 1/2$. Moreover, by Proposition \ref{prop:p3p5} these are in 1-to-1 correspondence with each other by replacing all $\psi_3$ by $\psi_8$ and conversely. It suffices therefore, to show that any $\psi$ from $ \mathcal{M}_{3,8}\smallsetminus \mathcal{M}_3 $ starting with $\psi_3$ fixes an  $s_{\alpha, 0}=0c_\alpha$.
This can be done with an argument similar to the one above. Let $\psi=\psi_3\psi_{i_2}\dots\psi_{i_m}$. When $m=2$, $\psi=\psi_3\psi_8$, and we know that
$\psi(0c_\alpha)=c_\alpha$, where $c_\alpha$ is the fixed point of $\psi_1^2$. Proceed by induction, using Proposition \ref{prop:psi-g-d}.
Now $\psi$ can be written as $\psi'\psi''$ with $\psi'=\psi_3\dots$ and $\psi''=\psi_3\dots$ \emph{unless} $\psi$ has the form $\psi_3\psi_8^{m}$, but
then we can use Lemma \ref{lem:131}.

To handle $s'_{\alpha, 0}=1c_\alpha$, we use the property that in general
$ s'_{\alpha, \rho}=E\,s_{1-\alpha, 1-\rho}$ (see \cite[Lemma 2.2.17]{Lothaire}).
This yields $$ s'_{\alpha, 0}=E\,s_{1-\alpha, 1}=E\,s_{1-\alpha, 0}.$$
Since in general $E(w)$ is a fixed point of $E\sigma E$ when $w$ is a fixed point of $\sigma$, we obtain that the $ s'_{\alpha, 0}$ are generated by the morphisms from the monoid $\mathcal{M}_{4,7}$, since $E\psi_3 E=\psi_7$ and $E\psi_8 E=\psi_4$. \hfill\qed


\medskip

\noindent {\bf Remark 1} There is an interesting coding from the morphisms starting with $\psi_3$ in $\mathcal{M}_{3,8}$ to $\mathcal{M}^0_{1,3}$. Let + be binary addition with 3 and 8: 3+3=3, 8+8=3, 3+8=8, 8+3=8. Add $i_2i_3\dots i_m 3$ to $3\,i_2\dots i_m$, and replace $8$ by 1. For example: $38...88+88...83=83...38\mapsto 13...31$.


\medskip

We display the first three levels of the binary tree $\mathcal{T}_{3,8}$, where the nodes are labeled with the morphisms given by Theorem \ref{th:diam0}.
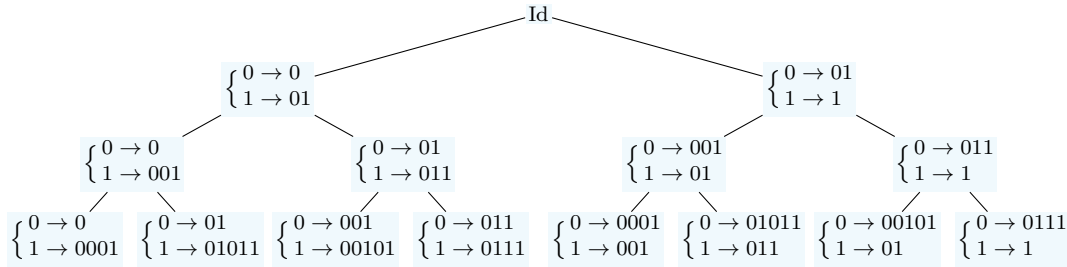
\begin{center}
\begin{tikzpicture}
[level distance=10mm,
every node/.style={fill=cyan!5,rectangle,inner sep=1pt},
level 1/.style={sibling distance=72mm,nodes={fill=cyan!5}},
level 2/.style={sibling distance=36mm,nodes={fill=cyan!5}},
level 3/.style={sibling distance=18mm,nodes={fill=cyan!5}}]
\node {\footnotesize  Id}            
child {node {\footnotesize $\mor{0}{01}$}            
 child {node {\footnotesize $\mor{0}{001}$}          
  child {node {\footnotesize $\mor{0}{0001}$}}       
  child {node {\footnotesize $\mor{01}{01011}$}}     
}
 child {node {\footnotesize $\mor{01}{011}$}         
  child {node {\footnotesize $\mor{001}{00101}$}}     
  child {node {\footnotesize $\mor{011}{0111}$}}    
}
}
child {node {\footnotesize $\mor{01}{1}$}            
 child {node {\footnotesize $\mor{001}{01}$}          
  child {node {\footnotesize $\mor{0001}{001}$}}       
  child {node {\footnotesize $\mor{01011}{011}$}}     
}
 child {node {\footnotesize $\mor{011}{1}$}         
   child {node {\footnotesize $\mor{00101}{01}$}}    
   child {node {\footnotesize $\mor{0111}{1}$}}  
}
};
\end{tikzpicture}
\end{center}

\medskip

\noindent {\bf Remark 2} Just as in Theorem \ref{th: hom}, each morphism generating an $s_{\alpha,0}$ occurs exactly once on the tree $\mathcal{T}_{3,8}$. This can be deduced from the fact that we have a coding between $\mathcal{M}_{3,8}$ and $\mathcal{M}^0_{1,3}$, but also because the monoid generated by the incidence matrices of $\psi_3$ and $\psi_8$ is free. Arnoux remarks that this can be derived in an elementary way (\cite[Lemma 6.5.14]{Arnoux}).

\section{Generating substitution invariant Sturmian words}\label{sec:gen}

There is a direct, more analytic way to find substitution invariant Sturmian words. We use an idea already considered by self-similarity expert Douglas Hofstadter in 1963 (\cite{Hof}).
To solve the fixed point equation $\psi(s_{\alpha, \rho})=s_{\alpha, \rho}$ for $\psi$, we can equivalently solve the fixed point equation
$$T_\psi(x,y)=(x,y)\quad {\rm for\;} 0<x, y<1,$$
where $T_\psi=T_{i_1}\dots T_{i_n}$ if $\psi=\psi_{i_1}\dots \psi_{i_n}$ with the  $i_k$ from some subset of $\{1,\dots,8\}$. Here the $T_i$ are two-dimensional fractional linear functions, such that
$$\psi_i(s_{\alpha, \rho})=s_{T_i(\alpha, \rho)}.$$
  Some $T_i$ are given by \cite[Lemma 2.2.18]{Lothaire}, and the others can be computed in a similar way. We have for example
\begin{equation*}T_1(x,y)=\left(\frac{1-x}{2-x},\frac{1-y}{2-x}\right), \quad
T_3(x,y)= \left(\frac{x}{1+x},\frac{y}{1+x}\right), \quad T_8(x,y)=\left(\frac{1}{2-x},\frac{y}{2-x}\right).
\end{equation*}

Note that both $T_3$ and $T_8$ leave the line $y=0$ invariant; this suggests the use of products of $\psi_3$ and $\psi_8$ to solve the equation $\psi(s_{\alpha, 0})=s_{\alpha,0}$, as we did in Theorem \ref{th:diam0}. We mention that the triple $T_1,T_3,T_8$ occurs in \cite{ItoYasutomi}, where they are used to connect two-dimensional continued fraction expansions to substitution invariant Sturmian words.

Solving the equation $T_\psi(x,y)=(x,y)$ is straightforward: there is a one-dimensional fractional linear function fixed point equation for $x$, which is quadratic, and then there is a linear equation for $y$, since one can show by induction on the number of $\psi_i$ in  $\psi$ that only $\pm y$ will occur in the second component of $T_\psi(x,y)$.

We mention that is some cases the equation is actually  $\psi(s_{\alpha, \rho})=s'_{T_\psi(\alpha, \rho)}$, but this can be dealt with by passing to the square of $\psi$, or by using Proposition \ref{prop:loz}.


\end{document}